

\documentstyle{amsppt}
\input epsf
\magnification=\magstep1
\hsize=6 truein
\hcorrection{.375in}
\vsize=8.5 truein
\parindent=20pt
\baselineskip=14pt
\TagsOnRight
\NoBlackBoxes
\footline{\hss\tenrm\folio\hss}

\centerline{Ergodicity of Mapping Class Group Actions}
\centerline{on Representation Varieties, I. Closed Surfaces}

\vskip1truein

\centerline{Doug Pickrell and Eugene Z. Xia}
\centerline{Department of Mathematics}
\centerline{University of Arizona}
\centerline{Tucson, AZ 85721}
\centerline{pickrell\@math.arizona.edu xia\@math.umass.edu}

\bigskip

\bigskip

\flushpar Abstract.  We prove that the mapping class 
group of a closed surface acts ergodically on connected 
components of the representation variety corresponding 
to a connected compact Lie group.  

\bigskip

\centerline{\S 0. Introduction.}

\bigskip

Throughout this paper we fix a connected compact Lie 
group $K$, and we let $dg$ denote the unique normalized 
Haar measure on $K$. 

Let $\Sigma$ denote a closed oriented surface with a fixed 
basepoint, and let $\Gamma_{\Sigma}=\pi_0(Aut(\Sigma ))$, the mapping class 
group.  The representation variety $Hom(\pi_1\Sigma ,K)$ has a 
canonical $\Gamma$$_{\Sigma}$-invariant measure class, the Lebesgue class 
of the set of nonsingular points, and it is well-known 
that this class is represented by a $\Gamma_{\Sigma}$-invariant measure 
([AB], [Go1]).  Our aim is to prove the following 

\proclaim{Theorem(0.1)} The group $\Gamma_{\Sigma}$ acts ergodically on 
the Lebesgue class of each connected component of 
$Hom(\pi_1\Sigma ,K)$.   
\endproclaim

Let $H^1(\Sigma ,K)$ denote the moduli space of representations, 
i.e. the quotient of $Hom(\pi_1\Sigma ,K)$ by the conjugation action 
of $K$. The following was proven by Goldman for $K$ locally 
isomorphic to $SU(2)\times T$, where $T$ is a torus, in [Go2].

\proclaim{Corollary(0.2)}The group $\Gamma_{\Sigma}$ acts ergodically on 
the Lebesgue class of each connected component of 
$H^1(\Sigma ,K)$. Thus the Lebesgue class of each component 
has an essentially unique $\Gamma_{\Sigma}$-invariant representative, 
the canonical symplectic volume element.
\endproclaim

As in most problems involving representation varieties, 
we use sewing techniques.  Our basic idea is to prove 
$(0.1)$ in initial cases, the one and two-holed tori, 
for $a.e.$ boundary condition; this is easier than dealing 
with every boundary condition, because we can use 
harmonic analysis for $K$.  When we sew, because we 
integrate, the measure-theoretic ambiguity is washed 
away, and we obtain the Theorem.  

The space $H^1(\Sigma ,K)$ is $\Gamma_{\Sigma}$-equivariantly filtered, the filter 
corresponding to a (conjugacy class of a) closed subgroup 
$K_1$ consisting of those representations which have image 
in $K_1$.  The filter corresponding to $K_1$ is isomorphic to 
$H^1(\Sigma ,K_1)$, modulo the action of the normalizer of $K_1$ inside 
$K$; it has a canonical $\Gamma_{\Sigma}$-invariant measure class, and each 
such class is represented by an invariant measure.  Our 
result implies that $\Gamma_{\Sigma}$ acts ergodically on connected 
components for classes corresponding to connected 
subgroups.  We have not addressed the case of 
nonconnected $K$; in this case $\Gamma_{\Sigma}$ does not generally act 
trivially on $\pi_0$ of the moduli space (consider a finite 
group), and it seems difficult to make an enlightening 
statement.  

The symplectic volume element on $H^1(\Sigma ,K)$ is the 
temperature $T\downarrow 0$ limit of the 2-dimensional Yang-Mills 
measure $d\nu_{\frac 1TYM_2}$ on the space of gauge equivalence 
classes of all (generalized) $K$-connections.  Elsewhere we 
will prove that this measure is ergodic with respect to 
its symmetry group, the group of area-preserving 
diffeomorphisms.  

\smallskip

\flushpar(0.2) Notation.  Given a Lie group $G$, we will 
always use left translation to trivialize the tangent 
bundle:  
$$TG\to G\times \frak g:v\vert_g\to (L_{g^{-1}})_{*}(v\vert_g).\tag 0.3$$
In this frame the commutator of two vector fields 
$x,y:G\to \frak g$ is given by
$$[x,y]\vert_g=dy(x)\vert_g-dx(y)\vert_g+[x(g),y(g)].\tag 0.4$$

The adjoint action $Ad:G\times \frak g\to \frak g$ is abbreviated to 
$Ad_g(x)=x^g$.  If $G$ acts on a space $X$, then $X^g$ denotes 
the fixed point set.  

\bigskip

\centerline{\S 1. Basic Notions and Sewing.}

\bigskip

For the purposes of this paper, we will need to 
consider a somewhat nonstandard kind of boundary 
condition for surfaces with boundary.  

Consider a connected compact oriented surface $\Sigma$ 
equipped with a basepoint, and the following additional 
structure:  each boundary component is linked to the 
basepoint by a path, and each boundary component $c$ is 
labelled with a $+$ or $-$, and a group element $k_c$ of $K$.  We 
interpret the sign to mean that the boundary component 
$c$ has an intrinsic orientation that agrees, or disagrees, 
with the induced orientation from $\Sigma$; the intrinsic 
orientation of the boundary component $c$ gives us a 
preferred generator for $\pi_1(c)\subset\pi_1(\Sigma )$, which, by slight 
abuse of notation, we will also denote by $c$.  We define 
$$Hom(\Sigma ,K)=\{g\in Hom(\pi_1(\Sigma ),K):g\vert_c=k_c,\forall 
c\in\pi_0(\partial\Sigma )\}.\tag 1.1$$
This space only depends upon the basepoint and paths to 
the boundary components up to homotopy.  The pure 
mapping class group $\Gamma_{\Sigma}$ does {\bf not} in general act on this 
space; only the subgroup generated by Dehn twists along 
curves which do not cross the paths from the basepoint 
to the boundary components will act; we denote this 
group by $\pi_0(Aut(\Sigma ))$.  

If $\Sigma$ is a closed surface, then we can form the quotient 
of $Hom(\Sigma ,K)$ by the global gauge action of $K$ by 
conjugation; the quotient is denoted by $H^1(\Sigma ,K)$.  In this 
case, $\pi_0(Aut(\Sigma ))=\Gamma_{\Sigma}$, the mapping class group.  

Let $s$ denote a separating oriented simple closed curve 
on $\Sigma$.  We suppose that the basepoint is on $s$, and we 
suppose also that $s$ does not cut any of the paths from 
the basepoint to the boundary components.  Let 
$\check{\Sigma}_k=\Sigma^{-}_k\bigsqcup\Sigma^{+}_k$ denote the disconnected object obtained by 
cutting along $s$ and attaching one $-$ and one $+$, and same 
group element $k$, to the new boundary components.  The 
Seifert-Van Kampen Theorem implies that the projection 
$p:\check{\Sigma}\to\Sigma$ induces an exact sequence 
$$0\to\langle s\rangle\to\pi_1(\Sigma^{-})*\pi_1(\Sigma^{+})@>{p_{
*}}>>\pi_1(\Sigma )\to 0.\tag 1.2$$
where $<s>$ denotes the normal subgroup generated by 
the element $s^{-1}*s$.  Hence we have the following 
elementary

\proclaim{Sewing Lemma(1.3)}Assume $\Sigma$ has a group 
element boundary condition.  Then there is a bijective 
correspondence 
$$Hom(\Sigma ,K)=\bigsqcup_{k\in K}Hom(\Sigma^{-}_k,K)\times Hom(
\Sigma^{+}_k,K),$$
where $g\leftrightarrow (g^{-},g^{+})$, $g^{-}(s)=g^{+}(s)=k$, $g^{
\pm}=g\vert_{\pi_1(\Sigma^{\pm})}$.  This 
correspondence is equivariant with respect to $\pi_0(Aut(\check{\Sigma }
))$, 
the group generated by Dehn twists along curves which 
cross neither $s$ nor the paths from basepoint to 
boundary components.  
\endproclaim

\bigskip

\centerline{\S 2. Initial Cases.}

\bigskip

The basic insight of this paper is that in all cases 
involving boundary, $\Gamma$-ergodicity is equivalent to 
$\Cal G$-ergodicity, where $\Cal G$ is a $continuous$ group of 
volume-preserving transformations, for $a.e.$ boundary 
condition.  The latter problem reduces to a calculation 
concerning infinitesimal transitivity.  

\bigskip

\flushpar\S 2.1.  The one-holed torus, with group element 
boundary condition.

\smallskip

\centerline{\epsfysize=3in            \epsffile{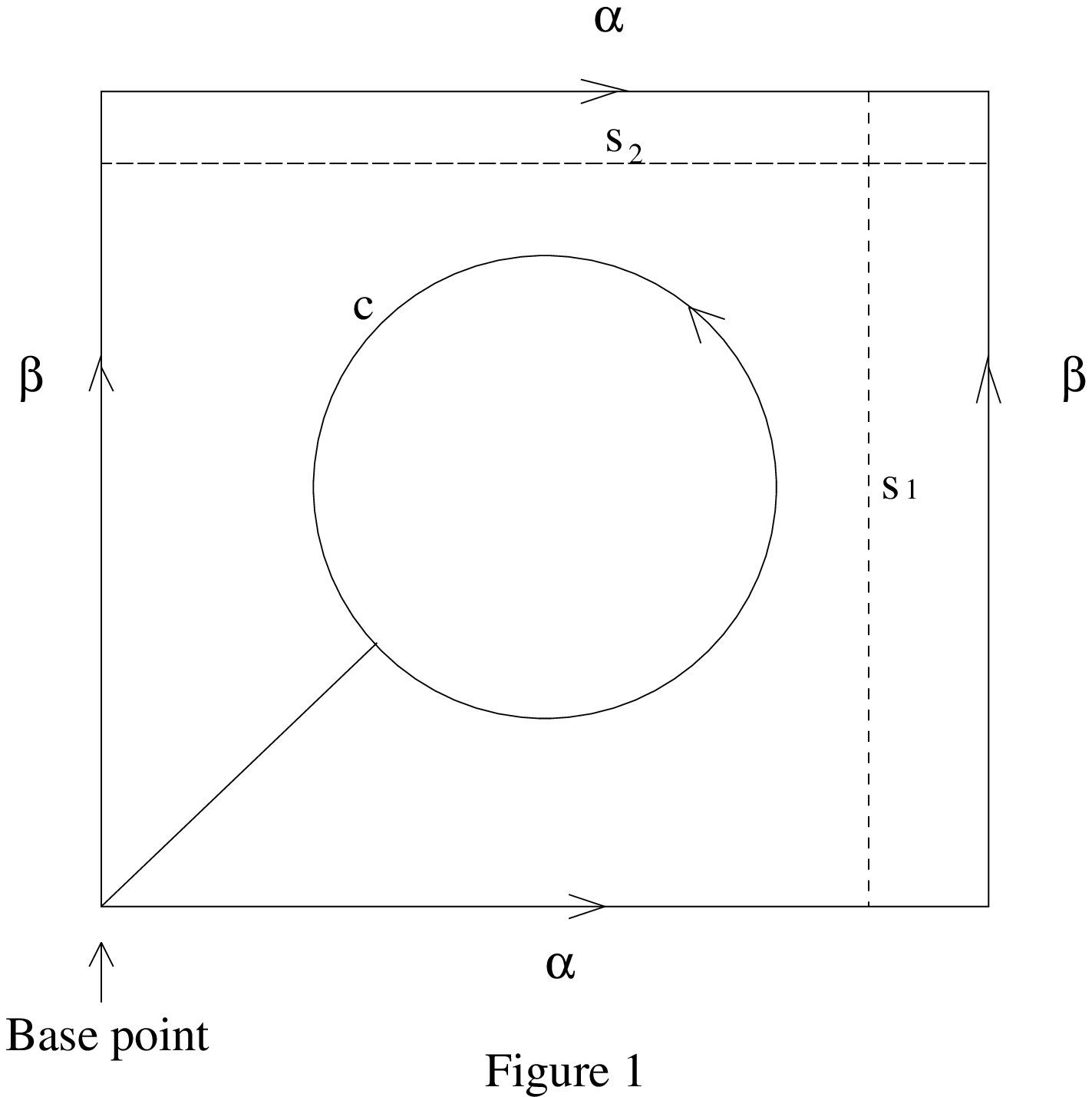}}

In this subsection we let $\Sigma$ denote the one-holed torus, 
with boundary component $c$, which we view as in Figure 
1.  

Given $k\in K$, we write $\Sigma_k$ to indicate that we impose the 
boundary condition $k$, so that $\Sigma_k$ is an object of the type 
considered in \S 1.  We have 
$$\matrix Hom(\pi_1\Sigma ,K)&\leftrightarrow&K\times K&(g\leftrightarrow 
(g_{\alpha},g_{\beta}))\\
&&\downarrow p\\
&&K'=[K,K]\endmatrix \tag 2.1.2$$
where $p$ is the commutator map, $p(g,h)=ghg^{-1}h^{-1}$.  With 
respect to this identification, the fibers of $p$ are 
precisely the representation spaces $Hom(\Sigma_k,K)$.  Define $\Gamma$ 
to be the group generated by the transformations 
$T_j:K\times K\to K\times K$ given by 
$$T_1(g,h)=(gh^{-1},h),\quad T_2(g,h)=(g,hg^{-1}).\tag 2.1.3$$
These transformations arise from twists along the 
curves $s_1$ and $s_2$ indicated in Figure 1; they are 
volume-preserving (with respect to Haar measure), hence 
naturally induce unitary transformations of $L^2(K\times K)$, 
they commute with conjugation by $K$, and they fix the 
map $p$.  The action of $\Gamma$ restricts to the action of 
$\pi_0(Aut(\Sigma_k))$ on the fiber $p^{-1}(k)=Hom(\Sigma_k,K)$.

In this subsection we will prove the following result, 
which is of independent interest.  

\proclaim{(2.1.4)Theorem} Suppose that $F\in L^2(K\times K)$ is 
$\Gamma$-invariant.  Then $F$ is $a.e.$ constant on components 
of $p^{-1}(k)$ for $a.e.$ $k$ $[d\rho ]$, where $d\rho =p_{*}(dg\times 
dh)$.  
\endproclaim

\flushpar Remarks (2.1.5).  (a) The measure $d\rho$ is in the
Lebesgue class of $K'$, and 
$$d\rho (k)=(\sum_{\mu}d^{-1}_{\mu}\chi_{\mu}(k))dk,\tag 2.1.6$$
where the sum is over all irreducible characters of $K'$, 
and $d_{\mu}=\chi_{\mu}(1)$.  To see this, first note that because $
p$ is 
a $conj(K)$-equivariant map, and $dg\times dh$ is conjugation 
invariant, $d\rho$ is conjugation invariant.  Secondly, if 
$f=\sum c_{\mu}\chi_{\mu}$ is a central function, then 
$$\int fd\rho =\int_{K'}\int_{K'}f(ghg^{-1}h^{-1})dgdh=\sum_{\mu}
c_{\mu}\int \{\int\chi_{\mu}(ghg^{-1}h^{-1})dh\}dg$$
$$=\sum c_{\mu}\int\frac {\vert\chi_{\mu}(g)\vert^2}{\chi_{\mu}(1
)}dg=\int (\sum c_{\nu}\chi_{\nu}(g))(\sum d_{\mu}^{-1}\chi_{\mu}
(g^{-1}))dg,$$
which heuristically explains (2.1.6) (The third equality 
uses the well-known integration formula 
$\int\chi (xhyh^{-1})dh=\chi (x)\chi (y)/\chi (1)$, which follows from 
observing that the left hand side is a central function 
for $(x,y)\in K\times K$, and computing the expansion in terms 
of characters for $K\times K$).   Because characters are 
orthogonal, 
$$\int\vert\sum d_{\mu}^{-1}\chi_{\mu}(g)\vert^2dg=\sum d_{\mu}^{
-2}\int\vert\chi_{\mu}(g)\vert^2dg=\sum d_{\mu}^{-1}.$$
The Weyl dimension formula implies that this sum is 
finite, provided $\frak k$ does not have $su(2)$ factors (see 
below); hence in most cases, the density in (2.1.6) 
represents an $L^2$ function on $K$.  In general, if we fix a 
maximal torus and positive Weyl chamber, so that we 
can parameterize the representations by dominant 
integral functionals $\mu$, then the Weyl character formula 
implies that for $g\in T$, 
$$\sum d_{\mu}^{-1}\chi_{\mu}(g)=\Delta (g)^{-1}\sum_{w\in W}(-1)^{
l(w)}\{e^{i\rho}\sum_{\mu}d_{\mu}^{-1}e^{i\mu}\}^w(g),$$
where $W$ is the Weyl group, $l(w)$ is the length of $w$, and 
$\rho$ is half the sum of the positive roots.  If we write $\mu$ 
in terms of the fundamental dominant integral functionals, 
$\mu =\sum n_j\mu_j$, then the Weyl dimension formula implies 
$$d_{\mu}=\prod_{\alpha >0}\frac {\langle\mu +\rho ,\alpha\rangle}{
\langle\rho ,\alpha\rangle}\sim\prod_jn_j^{\langle\mu_j,2\rho\rangle}
,$$
so that $\sum d_{\mu}^{-1}e^{i\mu}=\sum d_{\mu}^{-1}e^{i\sum n_j\theta_
j}$ always represents an $L^2$ 
function with respect to the Haar measure of $T$.  Since 
this function is the boundary values of a holomorphic 
function on $(\Bbb C_{\le 1})^r$, it cannot vanish on a set of positive 
measure.  This explains the meaning of the density 
(2.1.6), and shows that $d\rho$ is in the Lebesgue class (see 
Appendix B for a more direct proof).  

In the case of $K=SU(2)$, the density in (2.1.6), as a 
function of $diag(z,z^{-1})$, $z=e^{i\theta}$, is given by 
$$\frac 1{z-z^{-1}}\{z\sum_{d\ge 1}d^{-1}z-z^{-1}\sum_{d\ge 1}d^{
-1}z^{-d}\}=-Im\{\frac {e^{i\theta}ln(1-e^{i\theta})}{sin\theta}\}$$
$$=ln(2-2cos\theta )^{1/2}+cos\theta\frac {arg(1-cos\theta -isin\theta 
)}{sin\theta}.$$

\smallskip

(b) Theorem (2.1.4) gives an algebraic characterization of 
functions $F(g,h)$ which have the form $f(ghg^{-1}h^{-1})$, in 
situations where all the fibers $p^{-1}(k)$ are connected, e.g.  
for $K$ simply connected (see Appendix A).  It seems to 
be unknown whether there might be a reasonable 
characterization for more general groups (e.g.  finite 
groups).  

\smallskip

\proclaim{(2.1.7)Corollary of (2.1.4)} For $a.e.$ boundary 
condition $k\in K$, the action 
$$\pi_0(Aut(\Sigma_k))\times Hom(\Sigma_k,K)$$
is ergodic on the Lebesgue class of each component. 
\endproclaim

To prove (2.1.4), we need to be able to analyze the 
transformations in (2.1.3).  Let $T=T_2$ denote the unitary 
transformation on $L^2(K_1\times K_2)$ corresponding to the second 
of these transformations, where we have introduced 
copies $K_1$ and $K_2$ of $K$, for notational clarity.  Recall 
that the Peter-Weyl Theorem asserts that there is a 
$K\times K$-equivariant isomorphism 
$$\bigoplus_{\mu}\Cal L(V_{\mu})\to L^2(K):(L_{\mu})\to f,\quad f
(g)=\sum_{\mu}dimn(\mu )^{1/2}tr_{\mu}(L_{\mu}\pi_{\mu}(g^{-1})),\tag 2.1.8$$
where $\Cal L(V_{\mu})$ denotes the space of linear transformations 
of $V_{\mu}$, the sums are over all irreducible representations, 
and the linear action of $(g_l,g_r)\in K\times K$ on these 
respective spaces is given by 
$$L_{\mu}\to\pi_{\mu}(g_l)L_{\mu}\pi_{\mu}(g_r)^{-1}$$
$$f(g)\to f(g_l^{-1}gg_r).$$

\proclaim{(2.1.9)Lemma} Via the isomorphisms
$$L^2(K_1\times K_2)=L^2(K_1;L^2(K_2))=\bigoplus_{\mu}L^2(K_1;\Cal L
(V_{\mu})),$$
$T=diag(T_{\mu})$, where $T_{\mu}$ is the multiplication operator
$$T_{\mu}:L^2(K_1;\Cal L(V_{\mu}))\to L^2(K_1;\Cal L(V_{\mu})):F_{
\mu}(g)\to F_{\mu}(g)\pi_{\mu}(g)^{-1}.\tag 2.1.10$$
In particular
$$L^2(K_1;\Cal L(V_{\mu}))^{T_{\mu}}=\{F_{\mu}:F_{\mu}(g)\vert_{(
V_{\mu}^g)^{\perp}}=0,a.e.g\},\tag 2.1.11$$
and if $F_{\mu}$ is $T_{\mu}$-invariant, then (viewed now as a 
function of two variables)
$$F_{\mu}(g,h)=F_{\mu}(g,ha(g)^{-1}),\tag 2.1.12$$
for any measureable function $a:K\to K$ such that 
$[a(g),g]=1$ $a.e.$.
\endproclaim

\flushpar Proof of (2.1.9).  The formula for $T_{\mu}$ is a direct 
consequence of the Peter-Weyl theorem, and the other 
statements follow directly from the formula for $T_{\mu}$.//  

\smallskip

Note that
$$A=\{a:K\to K:[g,a(g)]=1,\forall g\}\tag 2.1.13$$
is an abelian subgroup of the gauge group $Map(K,K)$.  
We will assume that the maps in $A$ are smooth, unless 
noted otherwise. It is probably not the case that
$A$ is a Lie subgroup, because the family of projections 
onto the subalgebras $\frak k^g$, $g\in K$, is not smooth. Nonetheless 
we will refer to  
$$\frak a=\{x:K\to \frak k:Ad_g(x(g))=x(g),\forall g\}.\tag 2.1.14$$
as the Lie algebra of $A$, because it has the crucial 
property
$$exp(\frak a)\subset A.$$

The group $A$ acts on $K_1\times K_2$ in two ways, corresponding 
to the actions (2.1.3), by 
$$A_1(a):(g,h)\to (ga(h)^{-1},h),\quad A_2(a):(g,h)\to (g,ha(g)^{
-1}),\tag 2.1.15$$
respectively.  Note that the transformations $T_i^n$, $i=1,2$, 
correspond to $a(k)=k^n$.  Note also that the 
transformations (2.1.15) are volume-preserving.  

We can restate (2.1.9) as 

\proclaim{Lemma(2.1.16)}The $L^2$ function $F(g,h)$ is 
$T_j$-invariant if and only if $F$ is $A_j$-invariant, for $j=1,2$ 
(Here we can require the maps in $A$ to be $C^{\infty}$, $C^0$, or 
merely measureable - the basic result is insensitive to 
this requirement).  
\endproclaim

Let $\Cal G$ denote the closure of the group of 
volume-preserving transformations of $K\times K$ generated by 
$A_1$ and $A_2$, inside the Lie group of all volume-preserving 
diffeomorphisms of $K\times K$ (it will turn out that,  for our 
purposes, we could just as well consider the closure in 
the group of all volume-preserving transformations, in 
the natural strong operator topology).  It is unclear 
whether $\Cal G$ is a Lie group, but it is useful to think in 
these terms, as we will now see.  The Lie algebra 
actions corresponding to (2.1.15) are given by the vector 
fields on $K_1\times K_2$ 
$$dA_1(x)\vert_{g,h}=(-x(h),0),\quad dA_2(x)\vert_{g,h}=(0,-x(g))
,\tag 2.1.17$$
respectively, for $x\in \frak a$.  These actions do not necessarily 
commute.  

\proclaim{Definitions(2.1.18)} (a) $\frak g_0$ is the Lie algebra of 
vector fields on $K\times K$ given by
$$\frak g_0=\{(x\vert_h,y\vert_g):x,y\in \frak a\};$$

(b) $\frak g$ is the Lie algebra of vector fields on $K\times K$ 
generated by the family of Lie algebras 
$$\{Ad_{\sigma}\frak g_0:\sigma\in A_1\quad or\quad A_2\}.$$
\endproclaim

The bracket for $\frak g_0$ is given by 
$$[(x_1,y_1),(x_2,y_2)]\vert_{g,h}=(dx_2(y_1)\vert_h-dx_1(y_2)\vert_
h,dy_2(x_1)\vert_g-dy_1(x_2)\vert_g).\tag 2.1.19$$
(see (0.2)).

Heuristically $\frak g_0$ is the Lie algebra corresponding to the 
group generated by the identity components of $A_1$ and 
$A_2$, while heuristically $\frak g$ is the Lie algebra corresponding 
to $\Cal G$.  In practice we will think of $\frak g$ as an 
$Ad(\Gamma )$-invariant Lie algebra containing $\frak g_0$.  

\proclaim{Lemma(2.1.20)} Assuming we require maps to be 
$C^{\infty}$, we have $exp(\frak g)\subset \Cal G$.  
\endproclaim

\flushpar Proof of (2.1.20). Suppose that $\xi =(y,x)\in \frak g_
0$. Now 
$exp\{t(0,x)\}\in A_2$ and $exp\{t(y,0)\}\in A_1$, $\forall t$. Thus
$$exp(\xi )=\lim_{n\to\infty}(exp((y/n,0)exp(0,x/n))^n\in \Cal G,\tag 2.1.21$$
because Trotter's product formula is valid for vector 
fields on a compact manifold.  Therefore for any $\sigma\in A_j$, 
$$exp(Ad_{\sigma}(\xi ))=\sigma exp(\xi )\sigma^{-1}\in \Cal G.\tag 2.1.22$$
Using Trotter's product formula (and the analogue for 
brackets) in the same way, we see that for sums and 
brackets of such vector fields, we again exponentiate 
into $\Cal G$.//  

\smallskip

Our goal now is to show that the Lie algebra $\frak g$ is 
infinitesimally transitive off a set of codimension $>1$ 
along a generic fiber of the commutator map $p$.  We 
calculate that 
$$dp\vert_{g,h}:\frak k\oplus \frak k\to \frak k':(\xi ,\eta )\to
\xi^{hgh^{-1}}-\xi^{hg}+\eta^{hg}-\eta^h,$$
$$=(\xi^{h^{-1}}-\xi +\eta -\eta^{g^{-1}})^{hg}.\tag 2.1.23$$

\proclaim{Proposition(2.1.24)}For $g$ and $h$ in the 
complement of a set of codimension $>1$, the evaluation 
map 
$$eval\vert_{g,h}:\frak g\to ker(dp\vert_{(g,h)})$$
is surjective, where $eval$ is the evaluation map.
\endproclaim

\flushpar Proof of (2.1.24). If $K$ is abelian, then at all points 
$$eval\vert_g:\frak a\to \frak k^{Ad(g)}\tag 2.1.25$$
is surjective, and it follows from this that $\frak g$ is transitive.

So suppose that $K$ is nonabelian.  Recall the set of 
regular points, 
$$K^{reg}=\{k\in K:dimn(\frak k^g)=r\}\tag 2.1.26$$
where $r=rank(\frak k)$ is the minimal possible dimension of 
$\frak k^g$.  The singular set $K\setminus K^{reg}$ has codimension 3, because 
for a nonregular point $g$, $\frak k^g$ always contains a copy of 
$su(2)$, in addition to a maximal torus.  For a regular 
point $g\in K$, $eval\vert_g$ in (2.1.25) will be surjective (while 
the image shrinks at nonregular points, e.g.  
$eval\vert_1(\frak a)=\{0\}$); this follows from the real analyticity of 
the vector bundle $g\to \frak k^g$ over $K^{reg}$.  Therefore for 
$g,h\in K^{reg}$, 
$$eval\vert_{g,h}(\frak g_0)=\frak k^h\oplus \frak k^g\subset ker
(dp\vert_{g,h})\subset \frak k\oplus \frak k.\tag 2.1.27$$
This always fills out the central part of $\frak k$. For this 
reason, without loss of generality, we can henceforth 
assume that $\frak k$ is semisimple. 

The map $p$ is regular at all points $(g,h)$ such that 
$\frak k^g\cap \frak k^h=\{0\}$, by (2.1.23).  The abstract meaning of this 
condition is that the representation of $\pi_1(\Sigma )$ determined 
by $(g,h)$ is irreducible, in the intrinsic sense that the 
commutant of the image in $K$ is the center of $K$.  To 
understand this condition more concretely (from a point 
of view useful to us), suppose that $g$ and $h$ are regular.  
Write $h=exp(y)$, so that $\frak k^h=C_{\frak k}(y)$, the centralizer of $
y$.  
For $x\in \frak k^g$, 
$$[x,y]=\sum_{\alpha}\alpha (x)y_{\alpha}\tag 2.1.28$$
where the sum is over all roots $\alpha$ of $\frak k^g$, and $y_{
\alpha}$ denotes 
the component of $y$ in the root space of $\alpha$.  In order for 
$x\in \frak k^h$, we must have $\alpha (x)=0$, whenever $y_{\alpha}
\ne 0$.  The 
condition $y_{\alpha}=0$ is two independent real conditions, 
because the root space has one complex dimension.  Thus 
$\{(g,h):\frak k^g\cap \frak k^h\ne \{0\}\}$ has codimension at least 2.

We now know that the dimension of $ker(dp)$ is $dimn(\frak k)$ 
off a set of codimension 2. Let $proj_i$  denote projection 
onto the $i$th factor. The map $proj_1$ induces an exact 
sequence
$$0\to \{(0,\frak k^g)\}\to ker(dp)@>{proj_1}>>\{\xi\in \frak k:(
1-Ad(h^{-1}))\xi\in (\frak k^g)^{\perp}\}\to 0;$$
there is a similar sequence for $proj_2$.  The evaluation of 
$\frak g_0$ at $(g,h)\in K^{reg}\times K^{reg}$ fills out $ker(pr
oj_1)$$+$$ker(proj_2)$.  
Since 
$$\{\xi\in \frak k:(1-Ad(h^{-1}))\xi\in (\frak k^g)^{\perp}\}^{\perp}
=(1-Ad(h))\frak k^g,$$
to prove that $eval:\frak g\to ker(dp)$ is surjective at a regular 
point, it suffices to prove that 
$$(1-Ad(h))\frak k^g+proj_1(eval\vert_{g,h}(\frak g))=\frak k;\tag 2.1.29$$
there is a similar statement for $proj_2$.
 
Now $\frak g$ is $\Gamma$-invariant, hence  
$$\sum_{\Gamma}\gamma_{*}(eval\vert_{\gamma^{-1}(g,h)}(\frak g_0)
)\subset eval\vert_{g,h}(\frak g).\tag 2.1.30$$
In geometric terms, the sum is the $\Gamma$-invariant 
distribution generated by $\frak g_0$.  We will first consider only 
a small part of this sum, namely the $T_2$-invariant 
distribution generated by $\frak g_0$.  

Suppose that $(y\vert_h,x\vert_g)\in \frak g_0$.  We have
$$(T_2^n)_{*}(eval\vert_{T_2^{-n}(g,h)}(y,x))=$$
$$=(0,x\vert_g)+\frac d{dt}\vert_{t=0}(ge^{ty(hg^n)},hg^n(ge^{ty(
hg^n)})^{-n})$$
$$=\big\{\matrix (y(hg^n),x(g)-\sum_{k=1}^ny(hg^n)^{g^k}),\quad n
>0\\
(y(hg^n),x(g)+\sum_{k=n+1}^0y(hg^n)^{g^k}),\quad n<0\endmatrix .\tag 2.1.31$$
From this we see that 
$$\sum_{hg^n\in K^{reg}}\frak k^{hg^n}\subset proj_1(eval\vert_{g
,h}(\frak g)).\tag 2.1.32$$

Now for $(g,h)\in K^{reg}\times K^{reg}$, if $\{g^n\}$ is a dense subgroup 
of $T=exp(\frak t)$, then   
$$\sum_{\{n:hg^n\in K^{reg}\}}\frak k^{hg^n}=\sum_{\{x\in \frak t
:he^x\in K^{reg}\}}\frak k^{he^x}.\tag 2.1.33$$
For we clearly have $\subset$.  Conversely given $x\in \frak t$ such that 
$he^x$ is regular, we can find a sequence $\{n_j\}$ such that 
$g^{n_j}\to e^x$ as $j\to\infty$, hence $hg^{n_j}$ will be regular for $
j$ 
sufficiently large, and $\frak k^{hg^{n_j}}\to \frak k^{he^x}$, so that the opposite 
inclusion holds.

\proclaim{Lemma(2.1.34)} There is a set $X_2\subset K$ of 
codimension $\ge 2$ such that for $(g,h)\notin K\times X_2$, 
$$(1-Ad(h))\frak k^g+\sum_{\{x\in \frak k^g:he^x\in K^{reg}\}}\frak k^{
he^x}=\frak k.\tag 2.1.35$$
\endproclaim

\flushpar Proof of (2.1.34).  We write $h=exp(y)$.  We also 
write $\frak t=\frak k^g$.  Since $h$ is regular, there are open 
neighborhoods $\frak u$ and $U$ of $y$ and $h$, respectively, such 
that $exp:\frak u\to U$ is an isomorphism; let $log$ denote the 
inverse.  There is a Taylor series expansion of the form 
$$log(he^x)=\sum_{n\ge 0}c_n(h,x),\tag 2.1.36$$
where $c_n$ is homogeneous of degree $n$ in $x$.  If 
$\vert ad(y)\vert <\pi$, where $\vert\cdot\vert$ denotes the operator norm, then 
we can also expand each $c_n$, and the form of these 
expansions can be read off from the 
Baker-Campbell-Hausdorff formula, namely $c_0=y$, 
$$c_1=x+\frac 12[y,x]+\frac 1{12}[y,[y,x]]+..=x+O(\vert y\vert ),\tag 2.1.37$$
and for $n>1$,  
$$c_n=constant*ad(x)^n(y)+o(\vert y\vert ),\tag 2.1.38$$
as $\vert y\vert\downarrow 0$, where the constant depends only upon $
n$. We 
also have
$$(1-Ad(h))x=[x,y]+o(\vert y\vert )\quad as\quad\vert y\vert\downarrow 
0.\tag 2.1.39$$
We now claim that the sum in (2.1.35) equals
$$=(1-Ad(h))\frak t+\frak k^h+span\{log(he^x):x\in \frak t,he^x\in 
U\cap K^{reg}\}$$
$$=(1-Ad(h))\frak t+\frak k^h+span\{c_n(h,x):n\ge 0,x\in \frak t\}
.\tag 2.1.40$$
The first equality is immediate; the second follows 
from the fact that the span of the power series (2.1.36) 
will contain the span of the coefficients (replace $x$ by 
$sx$, note that $he^{sx}\in U\cap K^{reg}$ for small $s$, and 
differentiate with respect to $s$ at $s=0$).  

Now we first show that (2.1.35) holds for $h\in K^{reg}$ where 
$\vert y\vert$ is small. As in (2.1.28), $\frak k^h=C_{\frak k}(y
)$, and we can write
$$y=y_{\frak t}+\sum_{\alpha}y_{\alpha},\tag 2.1.41$$
relative to the root decomposition of $\frak k^{\Bbb C}$ with respect to 
$\frak t$.  If all the $y_{\alpha}\ne 0$, then together $\frak t$ and 
$\{ad(x)^n(y):x\in \frak t,n\ge 1\}$ will span $\frak k$.  But (2.1.37), (2.1.38) and 
(2.1.39) now imply that $\{c_n(h,x):n\ge 1,x\in \frak t\}$ and $(
1-Ad(h))\frak t$ 
will span $\frak k$, provided that $\vert y\vert$ is small.

Note that the condition $y_{\alpha}=0$ is linear, and of 
codimension 2, as we pointed out below (2.1.28).  Thus 
for $h$ in a subset of codimension $2$ in a neighborhood of 
$1$, the equation in (2.1.35) holds. It remains to do a 
similar analysis for a neighborhood of a point $e^{y_0}\ne 1$.  

Suppose that $h=e^{y_0+z}$, where $log(e^{y_0})=y_0$ and $z$ is 
small.  We have $c_0(h,x)=y_0+z$, 
$$c_1(h,x)=c_1(e^{y_0},x)+\frac 12[z,x]+o(\vert z\vert )=c_1(e^{y_
0},x)+O(\vert z\vert )\tag 2.1.42$$
$$(1-Ad(h))x=(1-Ad(e^{y_0}))x+[x,z]+o(\vert z\vert )\tag 2.1.43$$
$$c_n(h,x)=c_n(e^{y_0},x)+constant*ad(x)^n(z)+o(\vert z\vert )\tag 2.1.44$$
as $\vert z\vert\downarrow 0$, where $n>1$.  The derivative of $\frak k^
h$ (as it varies 
in the Grassmannian of subspaces, $Gr(r,\frak k)$) is a linear 
transformation $T(z):C(y_0)\to C(y_0)^{\perp}$.  If $T(z)(\xi_0)=
\xi_1$, then 
to first order in $s$, $exp(ad(y_0+sz))(\xi_0+s\xi_1)$$=$$\xi_0+s
\xi_1$, i.e.  
$[y_0,\xi_1]+[z,\xi_0]=0$; in terms of the root decomposition for 
$C(y_0)$, we have 
$$T(z)(\xi )=\sum_{\beta}\frac {\beta (\xi )}{\beta (y_0)}z_{\beta}
.\tag 2.1.45$$
Thus
$$\frak k^h=graph(T(z):C(y_0)\to C(y_0)^{\perp})+o(\vert z\vert )\tag 2.1.46$$
as $\vert z\vert\downarrow 0$.

Now consider the possibility that together $\frak k^{exp(y_0)}$, 
$(1-Ad(e^{y_0}))\frak t$, and $\{c_n(e^{y_0},x):n\ge 0,x\in \frak t
\}$ do not span $\frak k$.  
The argument proceeds initially as in the case 
$exp(y_0)=1$.  If $\forall\alpha$, $z_{\alpha}\ne 0$ (the components with respect 
to the root decomposition for $\frak t$), then 
$\{ad(x)^n(z):n\ge 1,x\in \frak t\}$ will span $\frak t^{\perp}$.  We now use 
(2.1.42)-(2.1.44).  For the variation of the span of $\frak k^h$, 
$(1-Ad(h))\frak t$, and $\{c_n(h,x):n\ge 0,x\in \frak t\}$ to be all of $
\frak k$, it is 
therefore sufficient for the natural map of an 
$r+1$-dimensional space to an $r$-dimensional space 
$$\Bbb Rz+graph(T(z))\to \frak k/\frak t^{\perp}\tag 2.1.47$$
to be surjective (note that the $z$ comes from the $c_0$ 
term; see the line preceding $(2.1.42)$).  Thus if (1) $z_{\alpha}
\ne 0$, 
$\forall\alpha$, and (2) (2.1.47) is surjective, then for $h$ 
corresponding to small $z$, (2.1.35) will hold.  

We have already remarked that the first condition is of 
codimension 2.  From the formula (2.1.45) for $T(z)$, we 
see that $z$ is generically independent of $graph(T(z))$, and 
$graph(T(z))$ is generically transverse to $\frak t^{\perp}$.  Therefore 
the condition (2) also has codimension 2, for $z$ in a small 
neighborhood of $1$.  This completes the proof.//  

\smallskip

We can now continue with the proof of (2.1.24).  Let $Q$ 
(for rational) denote the set of points of $K^{reg}$ with the 
property that $\{g^n\}$ is not dense in the torus $exp(\frak k^g)$.  

Now by (2.1.33) and (2.1.34) we know that $eval\vert_{(g,h)}$ maps 
onto $ker(dp)$ provided that $(g,h)$ is not in the set 
$(Q_1\times K_2)\cup (K_1\times X_2)$.  Now by considering the $T_
1$-invariant 
distribution generated by $\frak g_0$ (and $proj_2$), we can also 
conclude that $eval\vert_{(g,h)}$ maps onto $ker(dp)$ provided that 
$(g,h)$ is not in the set $(K_1\times Q_2)\cup (X_1\times K_2)$. 

Because the condition that $eval\vert_{(g,h)}$ maps onto $ker(dp)$ is 
a linear independence condition, involving real analytic 
vector fields, the set of points where $eval$ does not map 
onto $ker(dp)$ is generically real analytic, and by the 
proceeding paragraph of codimension $\ge 1$.  Since $X_1$ and 
$X_2$ have codimension $\ge 2$, the only portion of the singular 
set identified in the previous paragraph which could a 
priori support an object of codimension one is $Q_1\times Q_2$.  
But $Q$ has Hausdorff dimension $d-1$, where $d$ is the 
dimension of $K$, hence $Q_1\times Q_2$ has Hausdorff dimension 
$2d-2$.  So the singular set must have codimension at 
least $2$.  This completes the proof of (2.1.24).//  

\smallskip

\flushpar Proof of Theorem (2.1.4).  Suppose that 
$F\in L^2(K\times K)$ is $\Gamma$-invariant.  By (2.1.16) $F$ is $
\Cal G$-invariant.  
Now given a generic point where $\frak g$ is infinitesimally 
transitive along the fiber, the $\Cal G$-orbit of that point will 
be open in the fiber.  For a generic fiber, the 
complement of these open sets has codimension $>1$, by 
(2.1.24).  Hence for a generic fiber, the $\Cal G$-orbits necessarily 
coincide with the components of the fiber.  Thus an 
invariant $F$ is locally constant on connected components 
of $a.e.$ fiber.  // 

\bigskip

\flushpar\S 2.2. The $n$-holed torus, with group element 
boundary condition.

\smallskip

\centerline{\epsfysize=3in            \epsffile{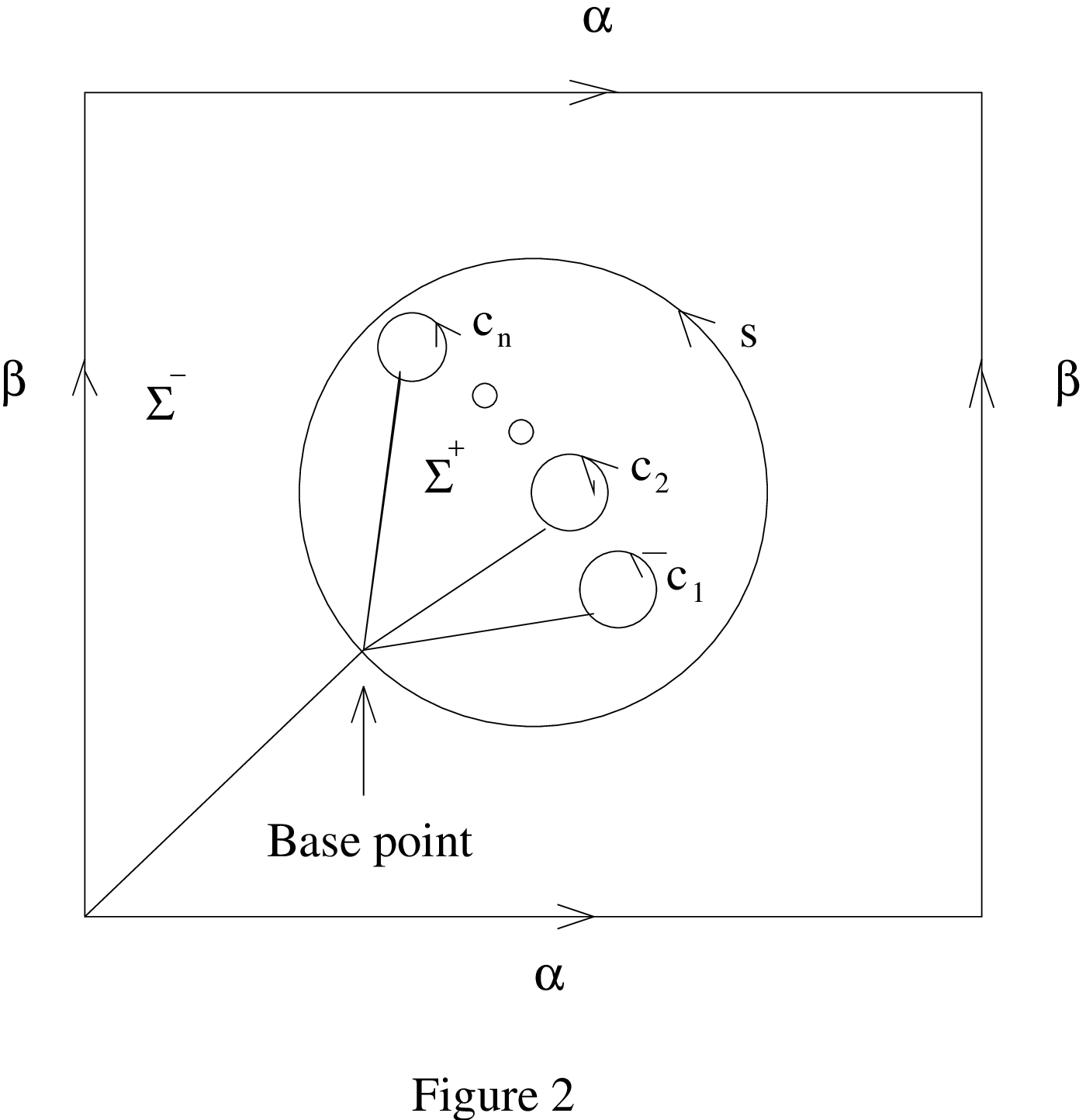}}

Let $\Sigma_{(k_1,..,k_n)}$ denote the $n$-holed torus with boundary 
components $c_1$$,..,c_n$, as in Figure 2 (where we ignore $s$ 
momentarily), with group element boundary condition 
$k_{c_j}=k_j$, $1\le j\le n$.  Let $\Sigma^{-}_k$ denote the one-holed torus 
with group element boundary condition considered in \S 2.1, 
which reappears in Figure 2 with boundary component $s$ 
(and where we have moved the basepoint from the 
vertex to $s$, which we can do without affecting the 
results of \S 2.1).  Let $\Sigma^{+}_{(k',k_1,..,k_n)}$ denote the $
n+1$-holed 
sphere with group element boundary condition pictured 
in Figure 2, where $k'$ is the labelling for $s$, and $k_j$ is the 
label for $c_j$; the corresponding $Hom$ space is empty 
unless $k'=\prod k_j$, in which case it is a point.  The Sewing 
Lemma (1.3) implies that we have a 
$\pi_0(Aut(\Sigma^{-}))$-equivariant bijection 
$$Hom(\Sigma_{(k_1,..,k_n)},K)\leftrightarrow Hom(\Sigma^{-}_k,K)
\times Hom(\Sigma^{+}_{(k,k_1,..,k_n)},K)\tag 2.2.1$$
where $k=\prod k_j$. Note that $\pi_0(Aut(\Sigma^{-}))$$=$$\pi_0(
Aut(\Sigma ))$.

Unfortunately this is not a situation where we can 
integrate over $k$, to obtain a result for every boundary 
condition, because $k$ is fixed by the $k_j$.  We need to vary 
one of the boundary conditions, say $k_n$.  We write $\Sigma_{(\vec {
k},\cdot )}$ 
for the object with boundary $k_{c_j}=k_j$, $1\le j<n$, where 
we allow $k_{c_n}$ to vary.  We then have a 
$\pi_0(Aut(\Sigma ))$-equivariant bijection 
$$Hom(\pi_1\Sigma^{-},K)\leftrightarrow Hom(\Sigma_{(\vec {k},\cdot 
)},K).\tag 2.2.2$$
An immediate consequence of Theorem (2.1.4) is the 
following 

\proclaim{Corollary(2.2.3)}For $a.e.$ $k_n$ $[d\rho ]$, the action 
$$\pi_0(Aut(\Sigma ))\times Hom(\Sigma_{(\vec {k},k_n)},K)$$
is ergodic on the Lebesgue class of each connected 
component.  
\endproclaim

\bigskip

\flushpar\S 3. Proof of Ergodicity.

\smallskip
 
Let $\Sigma$ denote the one-holed surface of genus $p$ with 
basepoint and link to the boundary $c$, as depicted in 
Figure 3 (ignore the paths $s$ and $\alpha$ at this point).  

\smallskip

\centerline{\epsfysize=4in            \epsffile{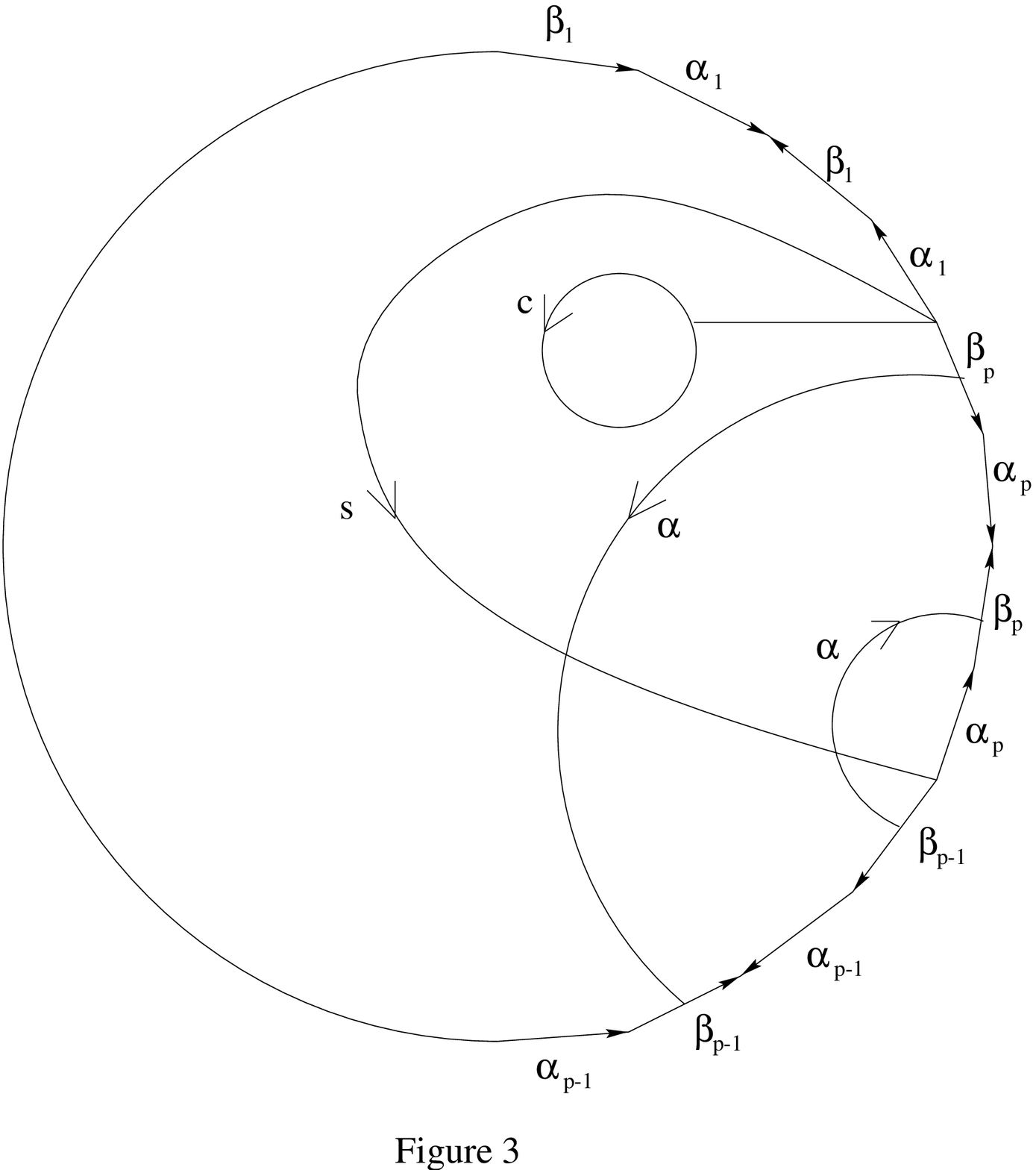}}

\proclaim{Theorem(3.1)}If the genus $p>1$, then for $\underline {
every}$ 
group element boundary condition $k\in K'$, the action 
$$\pi_0(Aut(\Sigma_k))\times Hom(\Sigma_k,K)\to Hom(\Sigma_k,K)$$
is ergodic with respect to the Lebesgue class of each 
connected component.  
\endproclaim

Note that Theorem (0.1) (when the $genus>1$) is the  
special case $k=1$ of $(3.1)$.  When the surface in 
(0.1) has $genus=1$, then (0.1) essentially reduces to the 
abelian case, and this is a standard application of 
Fourier series.  

The basic facts about the connectedness properties of 
$Hom(\Sigma ,K)$ which we will require are gathered in Appendix 
A, for the convenience of the reader.  In particular (A.3) 
asserts that $\pi_0(Aut(\Sigma ))$ acts on components, so that the 
statement of the Theorem makes sense.  

\smallskip

\flushpar Proof of (3.1).  Consider the decomposition of 
$Hom(\Sigma ,K)$ into connected components described in (A.3).  
If we prove (3.1) for all groups of the form $T\times K_1$, where 
$T$ is a torus and $K_1$ is simply connected, then we will be 
done.  So henceforth we assume that $K'$ is simply 
connected.  In this case all the representation spaces are 
connected, by (A.2).  

We have proven that $\pi_0(Aut(\Sigma^{-}))$ acts ergodically on 
$Hom(\Sigma^{-},K)$, for a one-holed torus $\Sigma^{-}$ as in \S 2.1, for $
a.e.$ 
group element boundary condition.  Similarly we have 
proven that $\pi_0(Aut(\Sigma^{+}))$ acts ergodically on $Hom(\Sigma^{
+},K)$ for 
a two-holed torus $\Sigma^{+}$ as in \S 2.2, for $a.e.$ boundary 
condition on one end, and for every boundary condition 
on the other end.  It therefore suffices to prove the 
following:  suppose that $s$ is a separating curve as in 
the Sewing Lemma (1.3), such that $\Sigma^{+}$ is a two-holed 
torus; if $\pi_0(Aut(\check{\Sigma }))$ acts ergodically on components of 
$Hom(\check{\Sigma },K)$, for $a.e.$ boundary condition on $s$, then the 
conclusion of (3.1) holds (see Figure 3; $\Sigma^{+}$ is to the 
reader's right of $s$).  

Let $k_c$ denote the fixed boundary condition for $\Sigma$.  The 
measure classes for possible boundary conditions on $s$ 
are the same for $\Sigma^{\pm}$, the Lebesgue class on $K'$ (see 
(2.1.5)).  Let $F$ denote a characteristic function on 
$Hom(\Sigma ,K)$.  If $F$ is $\pi_0(Aut(\check{\Sigma }))$-invariant, then by the 
Sewing Lemma (1.3) and our induction hypothesis, it 
follows that $F$ is constant along $a.e.$ fiber; hence $F$ is 
of the form $f(g\vert_s)$, where $f$ is a characteristic function 
on $K'$.  

Now suppose that $F$ is $\pi_0(Aut(\Sigma ))$-invariant.  As in 
Figure 3, we choose a Dehn twist $\sigma$ corresponding to a 
loop $\alpha$ that will cross the curve $s$, but will not cross 
the link to the boundary component, so that 
$\sigma\in\pi_0(Aut(\Sigma ))$.  Now the loop $\alpha$ does not pass through 
the basepoint.  There are two elemental ways in which 
we can use $s$ to link $\alpha$ to the basepoint; if we go from 
the basepoint in the negative direction along $s$ to $\alpha$, 
around $\alpha$, and return to the basepoint, then we denote 
this based loop by $\bar{\alpha}$; if we go from the basepoint in the 
positive direction along $s$ to $\alpha$, around $\alpha$, and return, 
then we denote the loop by $\underline {\alpha}$.  Using Figure 3 we 
compute that 
$$\bar{\alpha }=[\alpha_p,\beta_p]\alpha_p^{-1}\beta_{p-1}\alpha_{
p-1}\beta_{p-1}^{-1},\quad\underline {\alpha}=s\beta_{p-1}\alpha_{
p-1}\beta_{p-1}^{-1}\alpha_p[\alpha_p,\beta_p]^{-1}s^{-1},\tag 3.3$$
$$\sigma\circ\alpha_j=\alpha_j,\quad\sigma\circ\beta_{p-1}=\bar{\alpha}
\beta_{p-1},\quad\sigma\circ\beta_p=\beta_p\bar{\alpha}^{-1},\tag 3.4$$
$$\sigma\circ s=\underline {\alpha}s\bar{\alpha}^{-1}=(\prod_1^{p
-2}[\alpha_j,\beta_j])[\alpha_{p-1},\bar{\alpha}\beta_{p-1}]=\tag 3.5$$
$$s[\alpha_{p-1},\beta_{p-1}]^{-1}\alpha_{p-1}\alpha_p[\alpha_p,\beta_
p]^{-1}\alpha_{p-1}^{-1}[\alpha_{p-1},\beta_{p-1}][\alpha_p,\beta_
p]\alpha_p^{-1}.\tag 3.6$$

It is convenient to streamline our notation.  We put 
$$g_1=\prod_1^{p-2}[g_{\alpha_j},g_{\beta_j}],g=g_{\alpha_{p-1}},
h=g_{\beta_{p-1}},k=g_{\alpha_p},l=g_{\beta_p}.\tag 3.7$$
We have $(\sigma\cdot F)(g)=f(g_{\sigma\circ s})$. Hence the $\sigma$-invariance of $
F$ 
is equivalent to 
$$f(g_1[g,h])=f((g_1g)(kk_c^{-1})g_1[g,h](g_1g)^{-1}(kk_c^{-1})^{
-1}),\tag 3.8$$
for $a.e.$ $g_1,g,h,k,l$, subject to the constraint 
$g_1[g,h][k,l]=k_c$.

Define
$$\phi :\{(g_1,g,h,k,l):g_1[g,h][k,l]=k_c\}\to K'\times K'\tag 3.9$$
$$\phi :(g_1,g,h,k,l)\to (g_1[g,h],(g_1g)(kk_c^{-1})g_1[g,h](g_1g
)^{-1}(kk_c^{-1})^{-1});$$
(3.8) is equivalent to $\phi_1^{*}f=\phi_2^{*}f$, $a.e.$.  If this equality 
held at all points, then to prove that $f$ is constant, it 
would suffice to show that the relation defined by $Im(\phi )$ 
(or the equivalence relation generated by $Im(\phi )$) is 
transitive; since the equality holds in an $a.e.$ sense, we 
must consider the relation defined by the interior of 
$Im(\phi )$.  It is plausible that $\phi$ is surjective, but we can 
only prove the following weaker result.  

\proclaim{Lemma(3.10)} Let $pr_1:K'\times K'\to K':(m,n)\to m$. Then 
$K'\setminus pr_1(Interior(Im(\phi )))$ has codimension at least 2 in $
K'$.
\endproclaim

This Lemma implies that for each $m\in pr_1(Interior(Im(\phi )))$, 
we can find open sets $U_m$ and $V_m$ in $K'$ such that 
$m\in U_m$ and $U_m\times V_m\subset Im(\phi )$.  Since $\phi_1^{
*}f=\phi^{*}_2f$, $a.e.$, it 
follows that $f$ is constant on $U_m$, $a.e.$.  Since 
$pr_1(Interior(Im(\phi )))$ is connected, this constant must be 
the same for each $U_m$.  This implies that $f$ is constant 
$a.e.$.  Thus proving (3.10) will complete the proof of (3.1).  

\smallskip

\flushpar Proof of (3.10).  We can suppose that $p=2$, 
which amounts to setting $g_1=1$, and that $K=K'$.  For 
notational simplicity we will abbreviate $Ad(g)(\cdot )$ to $g(\cdot 
)$.  

To prove (3.10), we first claim that it suffices 
to show that the map 
$$\psi_m:\{[g,h]=m\}\times \{[k,l]=m^{-1}k_c\}\to K:(g,h;k,l)\to 
gkk_c^{-1}mg^{-1}(kk_c^{-1})^{-1}\tag 3.11$$
is regular at some smooth point, for each $m\in K\setminus Y$, 
where $Y$ has codimension 2.  For if $\psi_m$ is regular at the 
smooth point $(g,h;k,l)$, then $(g,h)$ is regular for the 
commutator map, which is the first factor of $\phi$.  Thus 
$Im(d\phi\vert_{(g,h,k,l)})$ spans both the vertical and horizontal 
directions, hence $(g,h;k,l)$ is regular for $\phi$.  

To specify $Y$, consider the commutator map $[,]:K\times K\to K$.  
This map is surjective, the fibers generically have 
dimension $d=dimn(K)$, and the exceptional fibers have 
dimension exceeding $d$ (e.g.  $[,]^{-1}(1)$ has dimension $d+r$, 
$r=rank(K)$).  Let $N$ denote the set of values $n\in K$ such 
that there exists $(g,h)$ with $[g,h]=n$ and (i) $g\in K^{reg}$ and 
(ii) $\frak k^g\cap \frak k^h=\{0\}$ (i.e.  $(g,h)$ is regular for $
[,]$).  By (B.5) of 
Appendix B, $K\setminus N$ has codimension at least 2 in $K$.  We 
set 
$$K\setminus Y=\{m:\quad m\in N\quad and\quad m^{-1}k_c\in N\}.\tag 3.12$$

The Zariski tangent space to $[,]^{-1}(m)$ at $(g,h)$ is given by 
$$T\vert_{(g,h)}=\{(x,y):x^{h^{-1}}-x+y-y^{g^{-1}}=0\}.\tag 3.13$$
Projection onto the $x$ factor induces the exact sequence
$$0\to \{(0,y):y\in \frak k^g\}\to T\vert_{(g,h)}\to \{x:(1-h^{-1}
)x\perp \frak k^g\}\to 0.\tag 3.14$$
If $(g,h)$ satisfies (i) and (ii), then $(g,h)$ is a smooth 
point, and the spaces in (3.14) have dimensions $r$, 
$d$, and $d-r$, respectively. Note that
$$\{x:(1-h^{-1})x\perp \frak k^g\}=((1-h)\frak k^g)^{\perp},\tag 3.15$$
and this space depends only upon $g$: since $hgh^{-1}=(g^{-1}m)^{
-1}$, 
$h$ is unique up to multiplication on the right by 
$\lambda\in C_K(g)$, and $\lambda$ acts trivially on $\frak k^g$.
 
Fix $m\in K\setminus Y$.  The derivative of the map $\psi_m$ is given by 
$$d\psi\vert_{(g,h;k,l)}:(x,y;z,w)\to (x^{m^{-1}(kk_c^{-1})^{-1}}
-x)^{kk_c^{-1}g}+(z^{k_c^{-1}gm^{-1}k_c}-z)^k$$
$$=(x^{(kk_c^{-1}m)^{-1}}-x+z^{m^{-1}k_c}-(z^{m^{-1}k_c})^{g^{-1}
m})^{kk_c^{-1}g}.\tag 3.16$$
Together with (3.14) and (3.15) this means that we must 
show that for suitable $g,h,k,l$, the sum of subspaces
$$(1-(kk_c^{-1}m)^{-1})((1-h)\frak k^g)^{\perp}+(1-g^{-1}m)m^{-1}
k_c((1-l)\frak k^k)^{\perp}\tag 3.17$$
is all of $\frak k$.   

Now to deal with (3.17), we need some control over 
solutions to the constraint equations $[g,h]=m$, 
$[k,l]=m^{-1}k_c$.  For this purpose, consider the equation 
$[g_1,h_1]=n$.  In (B.1) of Appendix B, we show that for any 
maximal torus $T$, there exists a solution $(g_1,h_1)$ with 
$g_1\in T$.  For $n\in N$, by dimensional considerations, $g_1$ is a 
finite multi-valued function of $T$ (see (a) of (B.6) for 
explicit equations).  Apply this to $n=m$.  Given $T$, we 
obtain solutions $[g,h]=m$.  We have $g^{-1}m$$=$$hg^{-1}h^{-1}\in 
hTh^{-1}$.  
Therefore we obtain a finite number of tori $hTh^{-1}$.  We 
claim that the multi-valued map 
$$\phi_m:\{Tori\}\to \{Tori\}:T\to hTh^{-1}\subset C_K(g^{-1}m)\tag 3.18$$
is surjective.  In a loose sense the inverse is $\phi_{m^{-1}}$, 
because $[hgh^{-1},h^{-1}]=m^{-1}$.  More precisely, given a torus 
$T_1$, apply the preceding to $m^{-1}$ and $T_1$ to obtain $(g_1,
h_1)$ 
with $[g_1,h_1]=m^{-1}$ and $g_1\in T_1$.  Define $(g,h)$ and $T$ so that 
$g_1=hgh^{-1}$, $h_1=h^{-1}$, $T=h_1T_1h_1^{-1}$.  Then $[g,h]=m$,  and 
$T_1=hTh^{-1}$.  This proves the claim.  

Similarly the multi-valued map 
$$\Phi_{m^{-1}k_c}:\{Tori\}\to \{Tori\}:T\to (kl)T(kl)^{-1}\subset 
C_K(kk_c^{-1}m),\tag 3.19$$
where $[k,l]=m^{-1}k_c$, $k\in T$, is surjective, and the inverse, 
again in a loose sense, is $\Phi_{k_c^{-1}m}$.  For given $T_1$ we can 
find $[k_1,l_1]=k_c^{-1}m$, $k_1\in T_1$.  Define $k=(k_1l_{}{}_1
)k_1(k_1l_1)^{-1}$, 
$l=(k_1l_1)k_1^{-1}(k_1l_1)^{-2}$, $T=(k_1l_1)T_1(k_1l_1)^{-1}$.  Then $
[k,l]=m^{-1}k_c$, 
$k\in T$, and $T_1=klT(kl)^{-1}$.  
 
Choose the pairs $(g,h)$ and $(k,l)$ such that $[g,h]=m$ and 
$[k,l]=m^{-1}k_c$, and such that both pairs satisfy (i) and (ii) 
above.  It may be necessary to consider perturbations of 
these pairs.  We will refer to perturbations which fix 
the constraints as admissible.  The conditions (i) and (ii) 
are stable under small admissible perturbations.  The 
space $[,]^{-1}(m)$ has dimension $d$, and for $g$ as above the 
possible $h$'s with $[g,h]=m$ form an $r$ dimensional set.  
Thus an admissible small perturbation of $(g,h)$ gives a 
smooth $d-r$ dimensional perturbation of $g$, the tangent 
space of which is described by (3.15).  The same 
comments apply to $(k,l)\in [,]^{-1}(m^{-1}k_c)$.  

Now consider the subspace represented by the first term 
in (3.17).  We first fix $g$ and $h$.  We claim that we can 
choose an arbitrarily small admissible perturbation of 
$(k,l)$ such that 
$$(1-(kk_c^{-1}m)^{-1})((1-h)\frak k^g)^{\perp}=$$
$$Im(1-(kk_c^{-1}m)^{-1})=(\frak k^{kk_c^{-1}m})^{\perp}.\tag 3.20$$
This will hold if we can arrange for $((1-h)\frak k^g)^{\perp}$ to 
intersect $ker(1-(kk_c^{-1}m)^{-1})$ trivially, i.e.
$$\frak k^{kk_c^{-1}m}\cap ((1-h)\frak k^g)^{\perp}=\{0\}.\tag 3.21$$
Because $k$ is regular, $kk_c^{-1}m=(kl)k(kl)^{-1}$ is regular.  Thus 
$\frak k^{kk_c^{-1}m}$ has dimension $r$, and $((1-h)k^g)^{\perp}$ has dimension $
d-r$ 
(because of condition (ii)).  By (B.1) and (3.19) we can find 
an arbitrarily small admissible perturbation of $(k,l)$ such 
that the intersection (3.21) will be zero.  

We now fix our choice of $(k,l)$.  We claim that we can 
find an arbitrarily small admissible perturbation of $(g,h)$ 
such that 
$$(1-g^{-1}m)m^{-1}k_c((1-l)\frak k^k)^{\perp}=(\frak k^{g^{-1}m}
)^{\perp}.\tag 3.23$$
The argument is essentially the same. It suffices to 
establish
$$\frak k^{g^{-1}m}\cap m^{-1}k_c((1-l)\frak k^k)^{\perp}=\{0\}.\tag 3.24$$
As before, $g^{-1}m$ is regular, because $g$ is regular. By (3.18) 
we can arrange this by an arbitrarily small admissible 
perturbation. 

We now have found $g,h,k,l$ such that the image of the 
subspace (3.17) equals 
$$(\frak k^{kk_c^{-1}m})^{\perp}+(\frak k^{g^{-1}m})^{\perp}=(\frak k^{
kk_c^{-1}m}\cap \frak k^{g^{-1}m})^{\perp},\tag 3.25$$
and this equality is stable under small admissible 
perturbations.  Again by (3.18) and (3.19) we can find a 
small admissible perturbation so that (3.25) will be all of 
$\frak k$. 

We have now proven that the map $\psi_m$ is regular at 
some smooth point for each $m\in K\setminus Y$, and as we observed 
at the beginning of the proof, this implies (3.10).  //

\bigskip

\centerline{Appendix A. Connectedness Properties.}

\bigskip

The following results can be deduced from [BR] (and 
perhaps elsewhere).  We record them here for the 
convenience of the reader.  

\proclaim{Lemma(A.1)} Suppose that $K$ is simply connected. 
Then $\{(g,h)\in K\times K:[g,h]=k\}$ is connected, $\forall k\in 
K$.
\endproclaim

If $C$ denotes the conjugacy class containing $k$, then there 
is a surjective map 
$$\{[g,h]=k\}\to \{[g,h]\in C\}/conj(K)$$
and the fibers are homogeneous spaces for $K$.  The 
fibers are connected because $K$ is connected, and by [BR] 
the moduli space corresponding to $C$ is connected because 
$K$ is simply connected.  This establishes (A.1) (It would 
clearly be desirable to give an elementary direct proof of 
this).  

\proclaim{Lemma(A.2)} Suppose that $K$ is simply connected.  
Suppose that $\Sigma$ is an object with group element boundary 
condition which is obtained by sewing one-holed tori to 
an $N$-holed sphere. Then $Hom(\Sigma ,K)$ is connected. 
\endproclaim

\flushpar Proof.  The space $Hom$ for an $N$-holed sphere is 
empty or a point.  When we sew, we obtain a connected 
object by (1.3).//  

\smallskip

Let $pr:\tilde {K}\to K$ denote the universal covering of $K$.

\proclaim{Proposition(A.3)} If $\Sigma$ is a one-holed surface 
with boundary condition $l\in K$, then we have the 
decomposition into connected components 
$$Hom(\Sigma_l,K)=\bigsqcup_{\tilde {l}\in\tilde {K}'\cap pr^{-1}
(l)}pr_{*}Hom(\Sigma_{\tilde {l}},\tilde {K}).$$
This decomposition is equivariant with respect to 
$\pi_0(Aut(\Sigma ))$.
\endproclaim

This follows from (A.2).

\bigskip

\centerline{Appendix B. Commutators.}

\bigskip

At several points of this paper, we used the fact that 
the commutator map $[,]:K\times K\to K'$ is surjective (and we 
presented an indirect proof of this in (a) of (2.1.5)).  
Here we discuss some refinements which we use in the 
proof of (3.10).  
  
\proclaim{Proposition(B.1)} Let $T$ denote a maximal torus 
in $K$. The map 
$$\psi :T\times K\to K':(\lambda ,h)\to [\lambda ,h]$$
is surjective.
\endproclaim
 
\flushpar Proof of (B.1).  To simplify the notation, we 
will write $K$ in place of $K'$; $d$ will denote the 
dimension, and $r$ the rank, of $\frak k$.

The derivative of $\psi$ at $(\lambda ,h)$ is given by 
$$\frak t\times \frak k\to \frak k:(x,y)\to (x^{\lambda h^{-1}}-x^{
\lambda}+y^{\lambda}-y)^h,\tag B.2$$
hence the image of the derivative at $(\lambda ,h)$ is
$$Ad(h)(Ad(\lambda )((1-Ad(h^{-1}))\frak t)+(\frak k^{\lambda})^{
\perp})$$
$$=Ad(h\lambda )((1-Ad(h^{-1}))\frak t+(\frak k^{\lambda})^{\perp}
).\tag B.3$$
We claim that the point $(\lambda ,h)$ is critical for $\psi$ if and 
only if (i) $\lambda\notin K^{reg}$ or (ii) $\frak k^h\cap \frak t
\ne \{0\}$.  To see this, 
suppose that $\lambda$ is regular and $\frak k^h\cap \frak t=\{0\}$.  Then $
(\frak k^{\lambda})^{\perp}=\frak t^{\perp}$ 
has dimension $d-r$ and $(1-Ad(h^{-1}))\frak t$ has dimension $r$.  If 
the intersection of these two spaces is nonempty, then 
there is $x\in \frak t$ such that $x^{h^{-1}}=x+y^{\perp}$, where $
y^{\perp}\in \frak t^{\perp}$ is 
not zero; but $x^{h^{-1}}$ and $x$ have the same length, so that 
$x\perp y^{\perp}$ implies $y^{\perp}=0$, which is a contradiction.  Thus 
the dimension of the space (B.3) is $d$, and this 
establishes our claim.  

 We can factor $\psi =\tilde{\psi}\circ p$, where 
$$T\times K@>{p}>>T\times K/T@>{\tilde{\psi}}>>K:(\lambda ,h)@>{p}>>
(\lambda ,hT)@>{\tilde{\psi}}>>[\lambda ,h],\tag B.4$$
so that at any regular point, $\tilde{\psi}$ will actually be a local 
diffeomorphism.  We claim that the set of critical 
values for $\tilde{\psi}$ has codimension at least two.  This will 
imply that $\tilde{\psi}$ is surjective, because a boundary for the 
image would necessarily have codimension one.  

Suppose that (i) holds, i.e.  $\lambda_0\notin T^{reg}$.  In this case, as 
we vary $h$, $\lambda_0h\lambda_0^{-1}h^{-1}$ will sweep out the $
\lambda_0$-translate of 
a nongeneric conjugacy class, which will have dimension 
$\le d-r-2$.  Thus the dimension of the set of critical 
values arising from condition (i) will be 
$\le r-1+d-r-2=d-3$.  

Now suppose that (ii) holds.  The condition $\frak k^{h_0}\cap \frak t
\ne 0$ 
has codimension at least 2 in $K$:  if $h_0=exp(X)$, where $X$ 
is regular, then we must have $X_{\alpha}=0$ for some root $\alpha$ 
of $\frak t$, where $X_{\alpha}$ denotes the $\alpha$-root space component of $
X$ 
(see the proof of (2.1.24), especially the paragraph 
containing (2.1.28)).  This is a $T$-invariant condition, 
hence the set of critical points corresponding to (ii) has 
codimension at least 2 in $K/T$.  It follows that the 
corresponding set of critical values has dimension 
$\le r+d-r-2$.  This completes the proof.//  

\smallskip

\proclaim{Corollary(B.5)} For the commutator map 
$[,]:K\times K\to K'$, the complement of the subset $N$ of $K'$ 
defined by 
$$\{n:\exists (g,h)\in [,]^{-1}(n)\quad s.t.\quad (i)\quad g\in K^{
reg},\quad (ii)\quad \frak k^g\cap \frak k^h=\{0\}\}$$
has codimension at least 2.
\endproclaim

\flushpar Proof.  Given a maximal torus $T$, each regular 
value for the map $\psi$ of (B.1) will belong to $N$.  In the 
proof of (B.1) we established that the complement of the 
set of regular values for $\psi$ has codimension at least 2.  
By varying $T$, we obtain (B.5).//  

\smallskip

\flushpar Remarks(B.6). (a) It is of interest to consider the 
more general question of whether, for given $g\in K$, the 
map 
$$\psi_g:T\times K'\to K':(\lambda ,h)\to [g\lambda ,h]\tag B.7$$
is surjective. This has a factorization
$$T\times K@>{p}>>D_g=\{(\lambda ,l)\in T\times K:g\lambda\sim l\}@>{
\tilde{\psi}_g}>>K@>{L_g}>>K$$
$$(\lambda ,h)@>{p}>>(\lambda ,hg\lambda h^{-1})=(\lambda ,l)@>{\tilde{
\psi}_g}>>\lambda l^{-1}=k@>{L_g}>>gk,\tag B.8$$
where $g\lambda\sim l$ means $g\lambda$ and $l$ are conjugate.  The map 
$\tilde{\psi}_g$ is the restriction to $D_g$ of the natural coset 
fibration 
$$T\times K\to (T\times K)/\Delta (T),\tag B.9$$
where $\Delta (T)$ is the diagonally embedded copy of $T$ in 
$T\times K$, and we identify $(T\times K)/\Delta (T)$ with $K$ by 
$(\lambda ,l)\Delta (T)\leftrightarrow\lambda l^{-1}$.  The map $
\tilde{\psi}_g$ is surjective if and only if 
for each $k\in K$, there exists $\lambda\in T$ such that $g\lambda
\sim l=k^{-1}\lambda$.  
This is equivalent to a system of $r$ polynomial equations 
$$\chi_i(g\lambda )=\chi_i(k^{-1}\lambda ),\quad i=1,..,r\tag B.10$$
for $r$ unknowns $\lambda_1,..,\lambda_r\in \Bbb T$, where $\chi_
i$ is the character 
corresponding to the $i^{th}$ fundamental irreducible 
representation, and $\lambda =\prod_1^r\lambda_i^{h_i}$, where the $
h_i$ are the 
coroots (e.g. for $SU(3)$, we have 2 equations
$$\sum_1^3A_i\lambda_i=0,\quad\bar {A}_1\lambda_2\lambda_3+\bar {
A}_2\lambda_1\lambda_3+\bar {A}_3\lambda_1\lambda_2=0,\tag B.11$$
for the $\lambda_i\in \Bbb T$, subject to the constraint $\prod\lambda_
i=1$, where 
$A_i=g_{ii}-(k^{-1})_{ii}$).  It is trivial to check that for $SU
(2)$, 
$\psi_g$ is always surjective, but this is not so for $SU(3)$.  
Thus in particular the equations (B.11) do not in general 
have solutions satisfying the reality condition $\vert\lambda_i\vert 
=1$; on 
the other hand (B.1) asserts that such solutions always 
exist for $g=1$.  

This suggests a number of questions, such as how does 
one describe the set of conjugacy classes which meet 
$gT$, when is $\psi_g$ surjective, and so on.  

(b)  Identify $SU(2)$ with $\Bbb H_1$, the group of 
unit quaternions, by $\left(\matrix a&b\\
-\bar {b}&\bar {a}\endmatrix \right)\leftrightarrow q=a-bj$, and take $
T=\Bbb T$.  
The conjugacy classes in $\Bbb H_1$ are obtained by fixing the 
real part of $q$.  Now fix $g=a-bj$.  The conjugacy 
classes which meet $g\Bbb T$ are indexed by $[-\vert a\vert ,\vert 
a\vert ]$.  We 
have 
$$D_g=\{(\lambda ,q)\in \Bbb T\times \Bbb H_1:\Bbb Re(q)=\Bbb Re(
a\lambda )\},$$
$$\tilde{\psi}_g:D_g\to \Bbb H_1:(\lambda ,q)\to\lambda\bar {q},$$
$$TD_g\vert_{\lambda ,q}=\{(is,q')\in i\Bbb R\times Im(\Bbb H):\Bbb R
e(a\lambda is)=\Bbb Re(qq')\},\tag B.12$$
$$d(\tilde{\psi}_g):TD_g\vert_{\lambda ,q}\to Im(\Bbb H):(is,q')\to 
q(is+\bar {q}')\bar {q},$$
$$D_{g,critical}=\{(\lambda ,a\lambda +\sqrt {1-\vert a\vert^2}zj
):\lambda ,z\in \Bbb T\}.$$
When $0<\vert a\vert <1$, the singular set is a 2-torus; at the 
extreme values $\vert a\vert =0,1$, the critical set degenerates to 
a circle.  The $SU(2)$ miracle is that in all cases, the set of 
critical values
$$\tilde{\psi}_g(D_{g,critical})=\{\bar {a}-\sqrt {1-\vert a\vert^
2}\lambda zj:\lambda ,z\in \Bbb T\}\tag B.13$$
is a circle.  One can easily visualize how $\tilde{\psi}_g$ covers $
\Bbb H_1$.  
  
The extreme case $\vert a\vert =0$, when there is just a single 
(totally geodesic) conjugacy class, corresponds to the 
condition that $g$ is a so-called principal element ([K]). 

\bigskip

\centerline{References}

\bigskip

[AB] M. Atiyah and R. Bott, The Yang-Mills equations 
over Riemann surfaces, Phil. Trans. R. Soc. Lond., A308, 
(1982), 523-615.

[BR] U. Bhosle and A. Ramanathan, Moduli spaces of 
principal bundles with parabolic structure over Riemann 
surfaces, Math. Z. 202 (1989), 161-180.  

[Go1] W. Goldman, The symplectic nature of fundamental 
groups of surfaces, Adv. Math. 54 (1984), 200-225.

[Go2] ----------, Ergodic theory on moduli spaces, Ann. 
Math. 146 (1997), 475-507. 

[K] B. Kostant, The principal three-dimensional subgroups 
and the Betti numbers of a complex simple Lie group, 
Amer. J. Math. 81 (1959), 973-1032.

\end